\documentclass[11pt]{article}
\usepackage{palatino}
\usepackage[letterpaper,left=1in,right=1in,top=0.75in,bottom=0.75in]{geometry}
\usepackage{natbib}
\usepackage{microtype}
\usepackage[hyperfootnotes=false]{hyperref}
\usepackage[font=small]{caption}
\usepackage[labelformat=empty, position=top]{subcaption}
\usepackage{graphicx}
\usepackage[export]{adjustbox}
\usepackage{picins}
\usepackage{titling}
\title{A simpler explanation of the reflective properties of conic sections, by following light rays along local isosceles paths}

\author{Rajeev D.\ S.\ Raizada}

\begin{document}

\vspace{-2.1in}

\maketitle

\parindent0cm
\parskip1.5ex
\renewcommand{\figurename}{Fig.}

\vspace{-0.5in}

\section*{Abstract}

Ellipses, parabolas and hyperbolas all have beautiful reflective properties. However, an intuitive explanation for {\em why} they have those properties has been lacking. There exist many mathematical proofs, but they tend to involve several analytical steps or geometrical constructions, making them unintuitive and hard to understand. Here, a simpler explanation is presented which only requires following the paths of light rays, and examining local paths that move from one point on a conic to a nearby one. First, a light-ray is followed as it runs along one of the legs of an isosceles triangle, and then reflects off a mirror that is parallel to the triangle's base. It bounces back along the path of the triangle's other leg. Next, a path is examined, moving from an arbitrary point on a conic section curve to a nearby point on the same curve. This path consists of two equal-length straight line steps, with each step following one of the constraints that defines the curve. For example, on an ellipse, defined by the sum of the distances to the two foci remaining constant, the path starts at a point on the ellipse, moves a distance $\delta$ directly away from one focus, then makes a second equal-length step directly towards the second focus. These two steps form the legs of precisely the sort of isosceles triangle described above, with its base running along the path of the curve. A light ray following the legs of that triangle gets reflected directly from one focus to the other. As the triangle shrinks towards zero, the reflection point converges onto the actual curve. Exactly the same argument also explains the reflections of parabolas and hyperbolas. Surprisingly, this explanation does not seem to have appeared previously in the long history of writings about conics. It is hoped that it will help to make the reflective properties of conic sections easier to understand and to explain.

\section*{Introduction: lots of proofs, but also lots of steps}

Amongst the best known and most important aspects of conic sections are their reflective properties. A parabolic mirror directs incoming parallel rays of light to its focus, as is exploited by every satellite dish, searchlight, and reflecting telescope. An ellipse sends any ray emanating from one focus towards the other focus, enabling medical devices called lithotriptors which focus ultrasound shockwaves to shatter kidney stones, and also giving rise to entertaining illustrations with elliptical billiard tables. A hyperbola reflects rays starting from one focus directly {\em away} from the second focus. This reflective property of a hyperbola turns out to be useful in an astronomical instrument called a Cassegrain telescope \citep[ch.6]{downs}, a noteworthy example of which is the Hubble Space Telescope \citep{nasa}.

There exist many ways of proving these properties, but they are not quite as straightforward as one might like. Perhaps for this reason, Algebra 2 textbooks which cover conic sections tend to describe the reflective properties without proving them \citep[e.g.][]{holliday,bittinger}. Proofs for the ellipse in particular often involve multiple steps. Perhaps the most commonly presented proof involves first invoking Fermat's principle that light travels along the path that takes the least time, i.e.\ the shortest possible one. Then, a tangent line is drawn at some point on the ellipse, and it is shown that the point of reflection along that tangent line that offers the shortest possible path is the one where it touches the curve of the ellipse itself, e.g.\ \citet[pp.6-8]{akopyan}. When restricted to reflection rather than also including refraction, this shortest-path rule is also referred to as Heron's principle \citep{fosterpedersen}.

Other proofs use trigonometry, e.g.\ the Law of Sines \citep[pp.29-30]{brannan}, or trigonometry together with calculating the derivative of the velocity vector of a point moving along the path of the curve \citep[pp.27-30]{glaeser}, or calculus combined with the trigonometric identity for the tan of the difference of two angles \citep[p.722]{stewart}.

Those sorts of proof, although of course sufficient to demonstrate that the reflective property holds true, can feel unsatisfying as explanations. In empirical studies, an explanation seeks not only to predict that some phenomenon will occur, but also to lay out what causes it to happen \citep{shmueli}. There is of course no causality in pure mathematics, but there is causality in the path of reflecting light rays. So, an intuitive explanation would show how, in the case of an ellipse, a light ray emanating from one focus gets caused to reflect back to the second focus.

Fermat's principle of least time is not easily intuitive, as it leaves somewhat mysterious the question of {\em how} the light manages to find the path that takes the least time. Although quantum mechanical explanations for this have been proposed \citep{feynman}, they are far from straightforward. (Perhaps also worth noting: for an elliptical mirror it seems that the principle would have difficulty explaining why the light ray should follow the direction along which it happens initially to be pointed. {\em Any} path from the first focus to the second one would be exactly as short as any other.) The proofs involving calculus and trigonometric identities do not suffer from that problem, but their multiple steps make it hard to form an intuitive understanding of how the premises eventually lead to the conclusion. 

To provide a causal explanation, the most direct approach would be to follow the path of a light ray and to see how it behaves when it meets the curve. For such an explanation to be intuitive, it would need to be simple and short. The approach proposed below attempts to satisfy those goals.

\section*{A simpler approach: bouncing light off the base of an isosceles triangle}

\piccaption{ 
\label{fig:isosc}}
\parpic[r]{\includegraphics[width=2.1in]{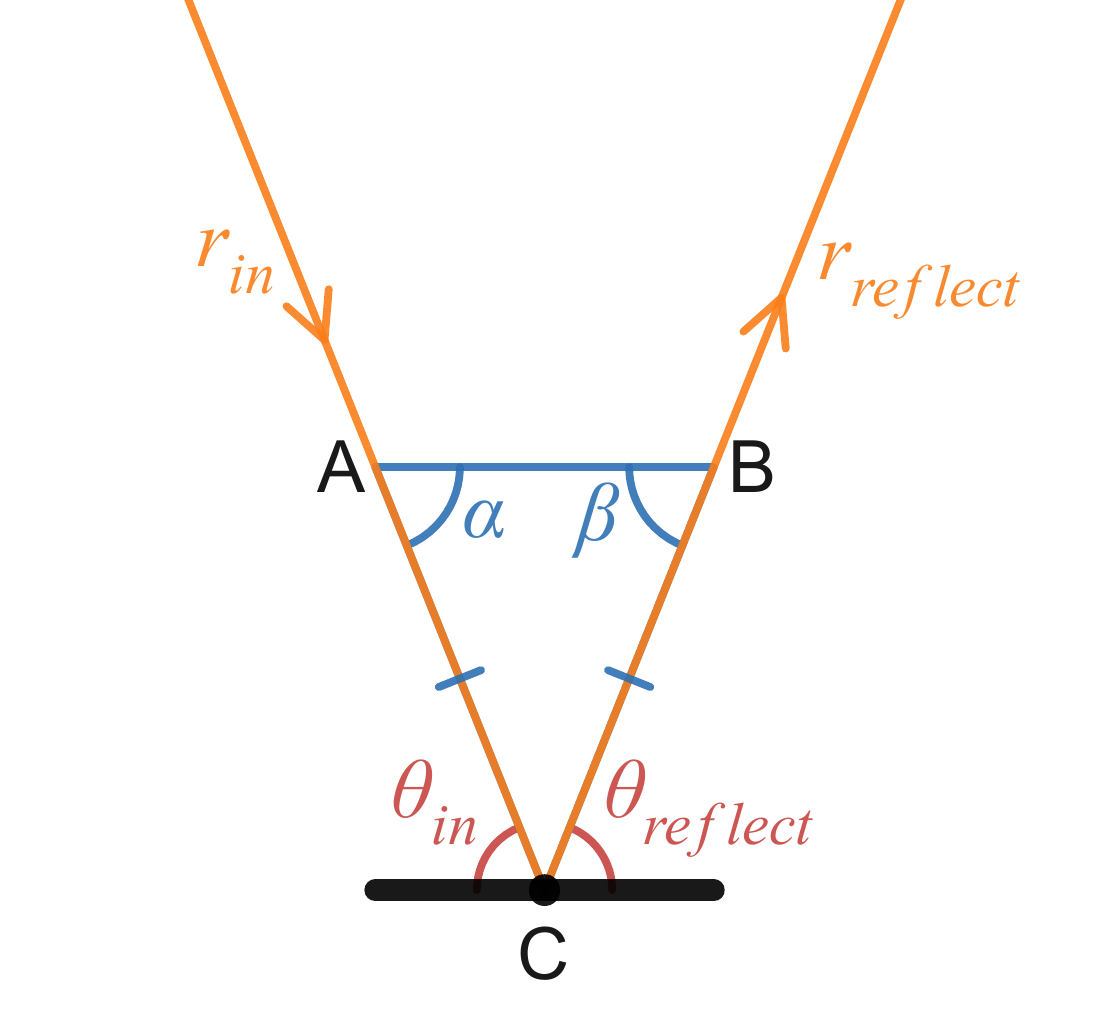}}

First, let us consider how light reflects off the base of an isosceles triangle, as shown in Fig.~\ref{fig:isosc}. In the triangle $ABC$, the legs $AC$ and $BC$ are of equal length. The base angles of this isosceles triangle, namely $\alpha$ and $\beta$, are therefore also equal. Let us reflect light off the thick line passing through the vertex $C$, parallel to the base $AB$ of the triangle. By the laws of reflection, the angle of incidence $\theta_{in}$ and the angle of reflection $\theta_{reflect}$ will be equal. Because the reflector is parallel to $AB$, the angles $\alpha$ and $\theta_{in}$ are alternate interior angles and hence are equal, as are $\beta$ and $\theta_{reflect}$. Thus, all four angles are equal to each other, and the reflected ray $r_{reflect}$ simply follows the path of the leg $BC$.

\section*{Moving along the curve of a conic by taking two equal-length steps}

\begin{figure}[!h]
  \centering
  \begin{subfigure}[t]{0.03\textwidth}
    \textbf{\Large A}
  \end{subfigure}
  \begin{subfigure}[t]{0.45\textwidth}
    \includegraphics[height=1.6in, valign=t]{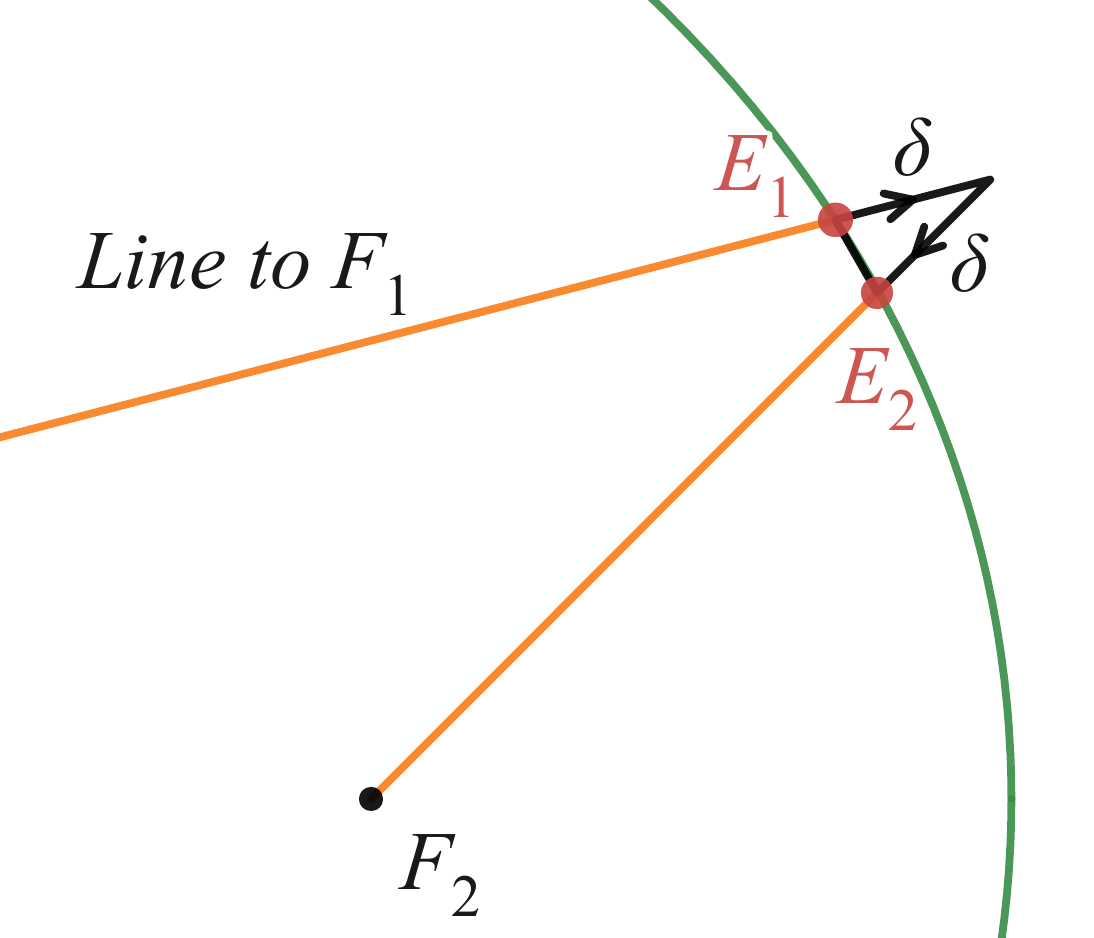}
  \end{subfigure}\hfill
  \begin{subfigure}[t]{0.03\textwidth}
    \textbf{\Large B}
  \end{subfigure}
  \begin{subfigure}[t]{0.45\textwidth}
    \includegraphics[height=1.6in, valign=t]{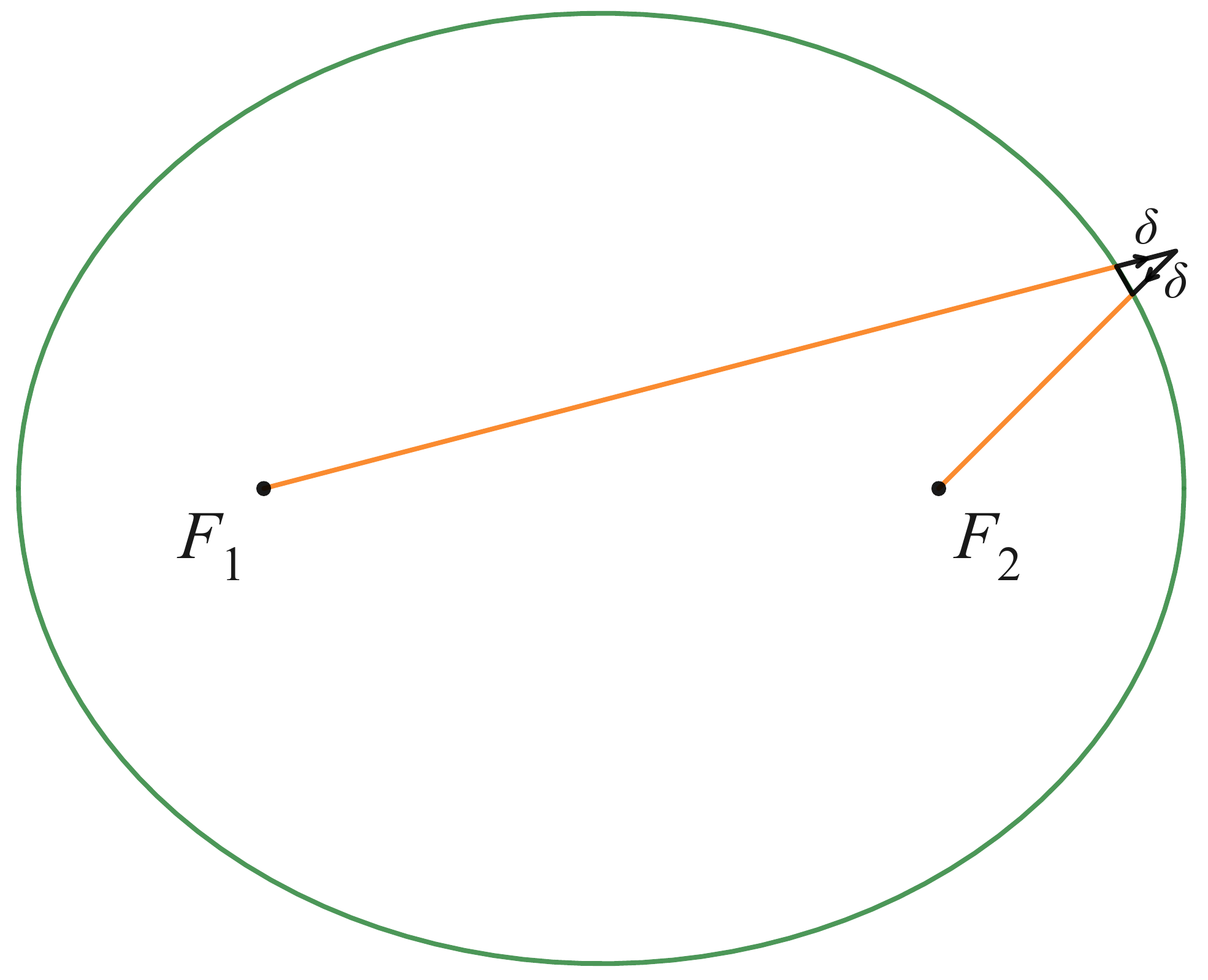}
  \end{subfigure}
  \caption{Drawing an isosceles triangle in two equal-sized steps, in order to move from one point to a neighbouring one along the curve of the ellipse. (A) A close-up view of the isosceles triangle formed by taking two steps of size $\delta$. (B) A broader view showing the entire ellipse.}
  \label{fig:ellipse}
\end{figure}

What does this seemingly pointlessly obvious fact about isosceles triangles have to do with the reflective properties of conics? Consider the question of how to move from one point on a conic section curve to a nearby point on the same curve. Figure~\ref{fig:ellipse} shows this for the case of an ellipse. Specifically, in Fig.~\ref{fig:ellipse}A, we move from point $E_1$ to point $E_2$, both lying on the ellipse. Recall that the ellipse is the set of all points such that the sum of the distances to the two foci is equal to some constant. So, if we start at $E_1$ and first move a distance $\delta$ away from the focus $F_1$, then we can keep this sum of distances constant and hence return to the ellipse by taking an equal-length step of $\delta$ towards the focus $F_2$, thereby arriving at the point $E_2$. 

The first step, moving directly away from $F_1$ by a distance $\delta$, added to the sum of the distances to the foci. The second step, of equal size but this time moving directly towards $F_2$, then subtracted an equal but opposite amount from that sum. So, the combined result of the two steps is to leave that sum unchanged, thus staying on the locus of the ellipse.

Crucially, these two equal-sized steps form an isosceles triangle, whose base $E_1 E_2$ runs along a small segment of the curve of the ellipse. This is exactly the sort of triangle discussed above and shown in Figure~\ref{fig:isosc}. A ray of light originating from focus $F_1$ and then bouncing off a line at the apex of the triangle that is parallel to its base will be reflected back along the second leg of the triangle, i.e.\ directly towards focus $F_2$.

When the triangle is of finite size, as shown in the figure, this reflection takes place very slightly outside of the ellipse, rather than directly on it. However, if we shrink the triangle ever smaller by letting the step-size $\delta$ tend towards zero, then the reflection tends towards happening exactly on the curve of the ellipse, and the base of the triangle tends towards being a perfect tangent. The infinitesimally small isosceles triangle receives light directly from the first focus and reflects it back exactly towards the second focus.

\newpage

\section*{The second step returns exactly to the path of the curve, as the triangle shrinks towards zero size}

\piccaption{The two equal steps' orthogonal projections onto each other are themselves of equal length. So, the two-step journey leaves unchanged the sum of the distances to the foci, thus returning to the ellipse. \label{fig:orthog}}
\parpic[r]{\includegraphics[width=4in]{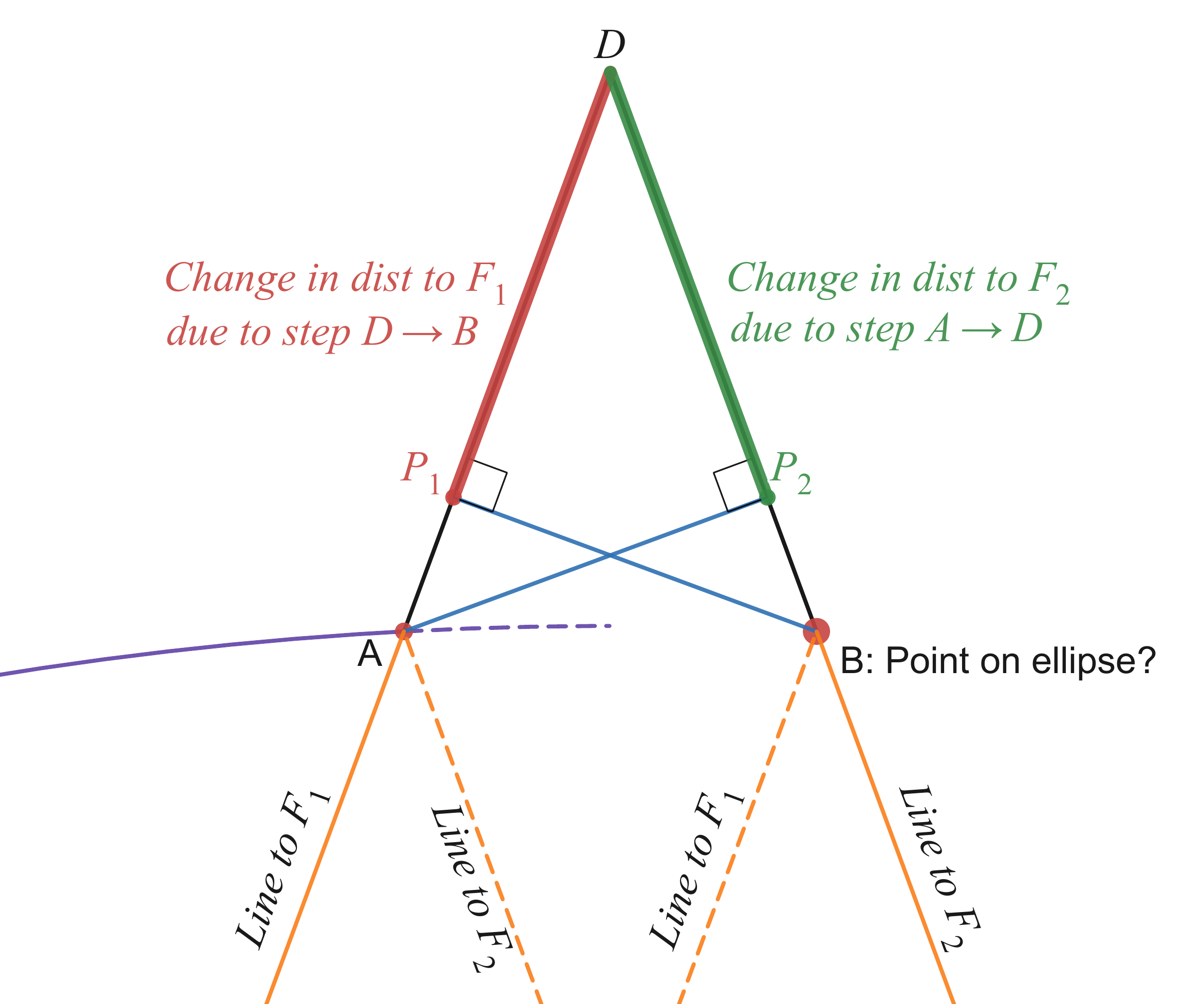}}

Before showing how the same argument also applies to parabolas and hyperbolas, it is worth pausing a moment to check that the second step of length $\delta$ truly does return to the locus of the ellipse. After all, although the first step does indeed move $\delta$ directly away from the focus $F_1$, it also adds some smaller but non-zero amount to the distance from the focus $F_2$. Similarly, the return step towards $F_2$ also reduces the distance to $F_1$.

Each such change is equal to the orthogonal projection of the $\delta$-length step onto the other leg of the isosceles triangle, as is shown in Figure~\ref{fig:orthog}. As can be seen from that figure, the symmetry of an isosceles triangle ensures that the two such projections are of equal length: the two steps $A \to D$ and $D \to B$ are both of length $\delta$, and they share the same projection angle $\angle ADB$. Hence, they are equal. So, the two $\delta$-length steps and the two equal-length projections collectively leave the sum of the distances to the two foci unchanged, thereby ensuring that the endpoint $B$ of the two steps does indeed lie on the ellipse, as required. 

Note that the above argument requires that the lines from $A$ to $F_2$ and from $B$ to $F_2$ must be parallel to each other, and similarly the lines from $A$ and $B$ to $F_1$. This is only approximately true for an isosceles triangle $ADB$ of finite size, but tends towards being exactly true as the size of that triangle shrinks towards zero compared to the distance to the foci. (One could construct a more formal argument that the lines become parallel, involving $\epsilon$ and $\delta$ and using the fact that $\tan \theta \to 0$ as $\theta \to 0$. However, given the aim of providing a simple explanation, and that it is clear that the lines would indeed tend towards becoming parallel, the additional complexity does not seem warranted.)

It is also worth asking whether this result might hold even if one of the legs of the isosceles triangle were {\em not} pointing towards a focus. To see that it would fail in such a circumstance, note that, for any given ellipse, its two foci are the {\em only} points for which the sum of their distances to the curve is constant. If the two $\delta$-length steps were to preserve the sum of distances to some other pair of points, then they would be the foci of a completely different ellipse.

\newpage

\section*{Applying exactly the same argument to parabolas and hyperbolas}

\piccaption{A parabola. \label{fig:parabola}}
\parpic[l]{\includegraphics[width=2in]{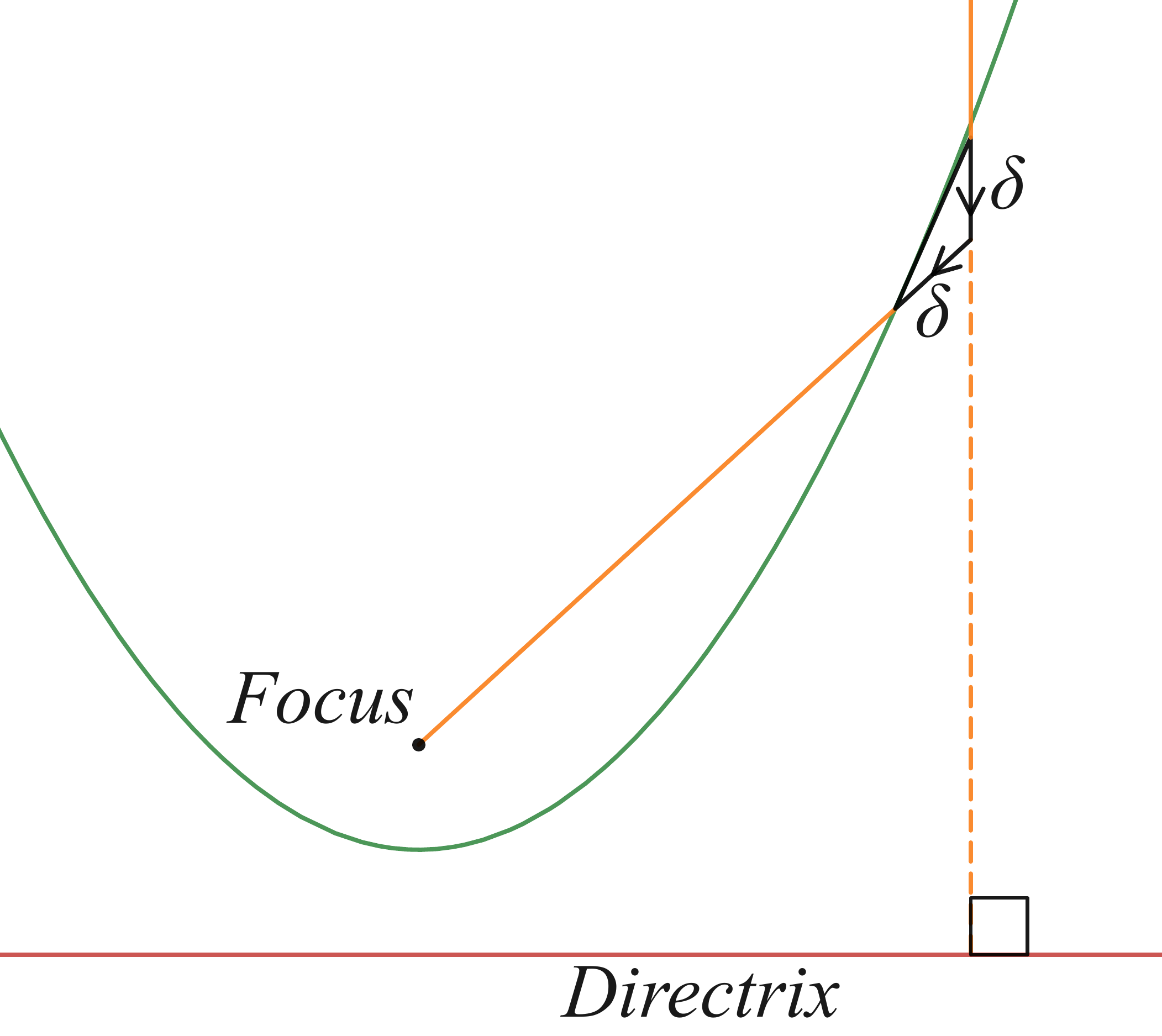}}

Exactly the same argument also applies to parabolas and hyperbolas, with only the directions of the two steps altering. Fig.~\ref{fig:parabola} shows two equal-length steps of $\delta$ being taken in order to move along the curve of a parabola. For points on this sort of curve, the distance to the directrix must remain equal to the distance to the focus. So, after taking a step of length $\delta$ straight towards the directrix, we can get back onto the parabola by taking an equal-length step straight towards the focus. 

The first step, moving directly towards the directrix, created an imbalance between the directrix- and focus-distances by adding $\delta$ to the former. The second step, also of length $\delta$ but now moving directly towards the focus, restored equality of the two distances by subtracting $\delta$ from the latter.

Here again, the resulting isosceles triangle shows why incoming light-rays that are perpendicular to the directrix will get reflected directly towards the focus. Also, as before, this triangle is only approximate when $\delta$ is finitely large, but it tends towards being exact as $\delta$ tends towards zero.

\piccaption{A hyperbola. \label{fig:hyperbola}}
\parpic[r]{\includegraphics[width=1.3in]{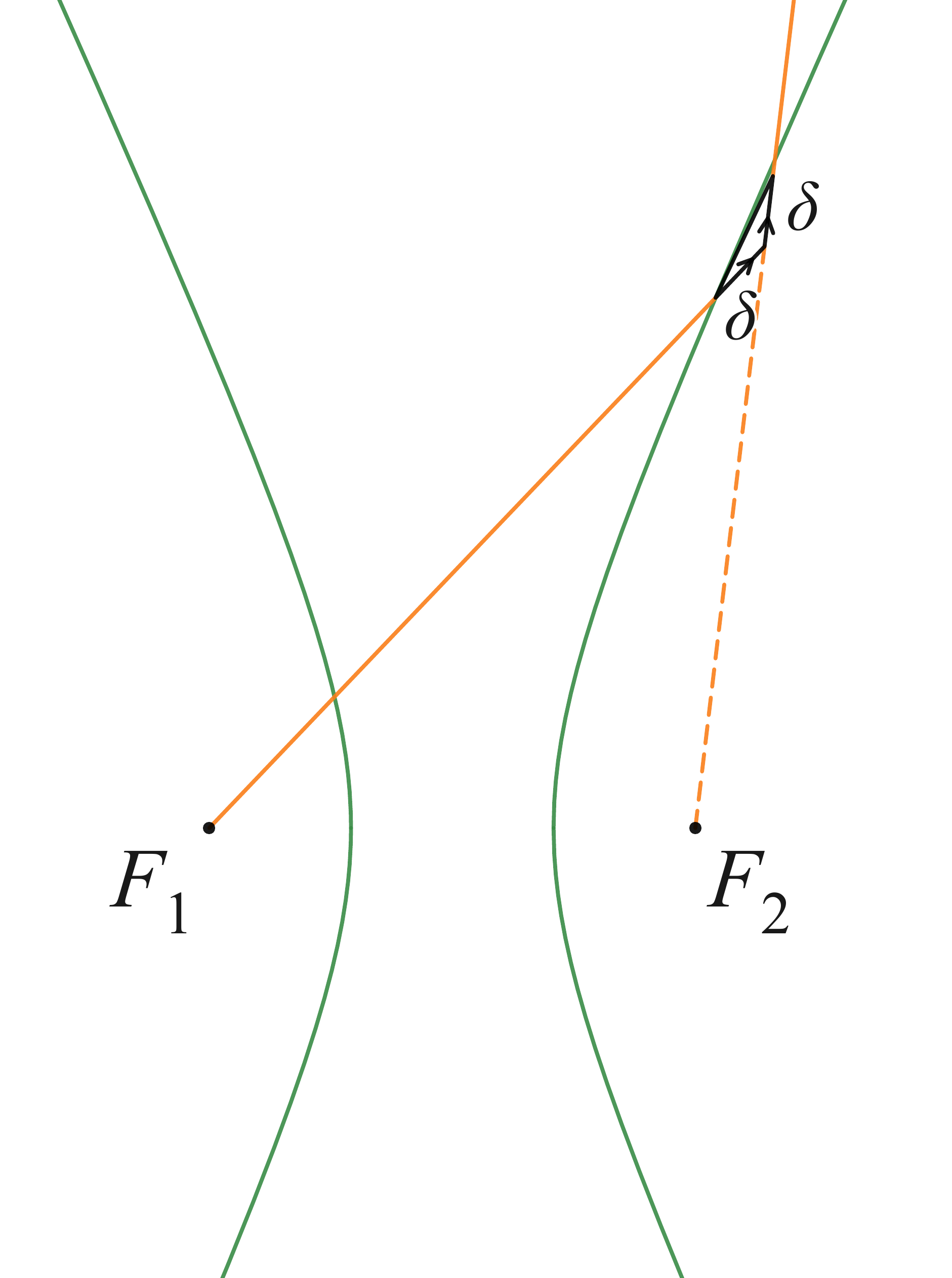}}

Figure~\ref{fig:hyperbola} similarly shows two equal-length steps of $\delta$ being taken in order to move along the curve of a hyperbola. In this case, the difference between the distances to the two foci must remain constant. So, after taking a step of length $\delta$ directly away from the first focus, $F_1$, we can get back onto the hyperbola by taking an equal-length step directly away from the second focus, $F_2$. For this type of curve, light emanating from one focus will get reflected directly away from the second focus, rather than towards it as was the case for the ellipse. As ever, the triangle shown in Fig.~\ref{fig:hyperbola} is only approximate, but the reflections become exact as $\delta$ tends towards zero.

\section*{Conclusion}

The argument presented here aims to provide a simple and intuitive explanation for the reflective properties of conics. It does not seem to exist in the current literature, either in the books on conics referenced above \citep{akopyan,downs,glaeser} or elsewhere. Whether previously published or not, the approach does not appear to be widely known. I hope that this simpler proof will be useful for teaching, and that it might help to make the reflective properties of conic sections easier to understand and to explain.

\bibliographystyle{apalike}
\bibliography{conics2}

\end{document}